\documentclass[a4paper,10pt]{article}
\usepackage[utf8x]{inputenc}
\usepackage{amsmath, amssymb, amsthm}
\usepackage{graphicx}
\usepackage{caption}
\usepackage{subcaption}
\usepackage{diagbox}
\usepackage{array}
\usepackage{tikz}
\usepackage{tabu}
\usepackage{enumerate}
\usepackage{float}
\usepackage{filecontents}
\restylefloat{table}
\usepackage[colorlinks,
  linkcolor=blue,
  citecolor=blue,
  filecolor=blue,
  menucolor=blue,
  urlcolor=blue,
  pdfnewwindow=true,
  pdftitle ={Maximal\ polytopes\ contained\ in\ a\ polytope},
  pdfauthor= {Moritz\ Firsching}, pdfsubject={Maximal\ polytopes\ contained\ in\ a\ polytope}]{hyperref}

\theoremstyle{definition}
\newtheorem{problem}{Problem}
\newtheorem{assumption}{Assumption}

\author{Moritz Firsching\thanks{Supported by DFG within the research training group ``Methods for Discrete Structures'' (GRK1408)}\\
\small Institut für Mathematik \\[-0.8ex]
\small FU Berlin\\[-0.8ex]
\small Arnimallee 2\\[-0.8ex]
\small 14195 Berlin\\ [-0.8ex]
\small Germany\\
\small \href{mailto:firsching@math.fu-berlin.de}{firsching@math.fu-berlin.de}
}
\title{Computing maximal copies of polytopes contained in a polytope}

\begin{document}
\maketitle

\begin{abstract}Kepler (1619) and Croft (1980) have considered largest homothetic copies of one regular polytope contained in another regular polytope. 
For arbitrary pairs of polytopes we propose to model this as a quadratically constrained optimization problem. 
These problems can then be solved numerically; in case the optimal solutions  are algebraic, exact optima can be recovered by solving systems of equations to very high precision and then using integer relation algorithms.
Based on this approach, we complete Croft's solution to the problem concerning maximal inclusions of regular three-dimensional polyhedra by describing inclusions for the six remaining cases.
\end{abstract}

\section{Introduction}
Given two polytopes $P$ and $Q$, we can ask: What is a polytope $P'$ of largest volume such that $P'$ is similar to $P$ and contained in~$Q$. 
By ``similar'' we understand that $P'$ can be transformed into $P$ by a dilation and rigid motions.
Instead of ``largest volume'' we might as well ask for a polytope that maximizes the dilation factor between $P$ and~$P'$. 
An equivalent question asks for the smallest polytope $Q'$, which is similar to $Q$ and contains~$P$. 

The earliest work on this topic might already be found in Kepler's work, \cite[libri V, caput I, p.~181]{K}. 
One finds descriptions of the largest regular tetrahedron included in a cube and of the largest cube included in a regular dodecahedron, although no claim on maximality is made.

A  substantial contribution is made by Croft, \cite{C80}. Here the case where $P$ and $Q$ are three-dimensional is considered. 
He notes that apart from exceptional cases \emph{local} maxima must be immobile and therefore satisfy $7$ linear constraints, see \cite[Theorem, p.~279]{C80}. 
Using this information he calculates all local maxima and obtains \emph{global} maximal configurations, see \cite[p.~283--295]{C80}. 
Letting $P$ and $Q$ range over the platonic solids, Croft gives a complete answer for $14$ out of the $20$ non-trivial cases. 
This is the problem described by the same author, Falconer and Guy as Problem B3 in \cite[p.~52]{CFG91}; see below for an answer for the remaining six cases. 

Containment problems for (simple) polygons are discussed for example in \cite{C83} and \cite{AAS98}, and some algorithms are given.
Taking $P$ to be a regular $n$-gon and $Q$ to be a regular $m$-gon, the size of the largest copy of $P$ inside $Q$ is known if and only if $n$ and $m$ share a common prime factor. 
If they are coprime only conjectural results are known; see the article by Dilworth and Mane, \cite{DM10}.

More general containment problems are studied by Gritzmann and Klee, \cite{GK94}. They also allow other groups than the group of similarities act on the polytopes. 
Gritzmann and Klee state the problem where the group acting is the group of similarities, \cite[p.~143]{GK94}, but do not discuss a computational approach.

The related problem of finding a largest, not necessarily regular, $j$-simplices in $k$-cubes is related to Hardamard matrices and discussed in \cite{HKL96}. In some cases the maximizer is indeed a \emph{regular} simplex, see \cite{MRT09} for details.

A short summary of the results of this paper by the author has been posted on mathoverflow, \cite{F13}.
\vspace*{.5cm}

In Section~\ref{methods} we present a method for finding solutions for this problem in general. 
In the last section we apply this method to some special cases and thereby offer a solution to Problem B3 in \cite[p.~52]{CFG91}.

\section{Methods}\label{methods}

\subsection{Setting up the optimization problem}
Let $P$ and $Q$ be polytopes, let $p$ be the dimension of $P$ and $q$ be the dimension of~$Q$. 
We assume $q\geq p$; otherwise it is not quite clear what it means that $P$ is included in~$Q$.
Let $H_1,\dots H_m$ be the defining half spaces for $Q$, such that 
\[Q=\bigcap_{k=1}^mH_k\]
and $w_1,\dots,w_n$ denote the vertices of~$P$.
We formulate the problem of finding the largest polytope $P'$ such that $P'$ is contained in $Q$ and similar to $P$ as a quadratic maximization problem.

\begin{problem}\label{prob}
{\centering\fbox{%
\begin{minipage}{.77\textwidth}
\begin{description}
 \item[Input data:] \[\text{halfspaces }H_1,\dots,H_m\text{ of }Q\text{, vertices }w_1\dots w_n\text{ of }P.\]%
 \item[Variables:] \[s\text{ and }v_{ij}\text{ for }1\leq i\leq n, 1\leq j\leq q\]
 \item[Objective function:] \[\text{maximize } s\]
 \item[Linear constraints:] \[(v_{i1},\dots,v_{iq})\in H_k \text{ for }1\leq i \leq n, 1\leq k\leq m\]
 \item[Quadratic constraints:] \[\sum_{l=1}^q(v_{il}-v_{kl})^2=s||w_i-w_j||^2_2\text{ for }1\leq i<j\leq n\]
\end{description}
\end{minipage}}}
\end{problem}
\vspace*{.1cm}
\noindent In this formulation the variable $s$ can be thought of as the square of the dilation factor between $P$ and~$P'$. The other variables are supposed to be the coordinates of the vertices of~$P'$. 
The linear constraints consist of $nm$ weak inequalities and make sure that $P'\subset Q$. 
The quadratic constraints assert that the distances between vertices of $P'$ agree with those of $P$ up to a dilation factor $\sqrt{s}$, which is the same for all pairs of vertices. 
Hence the quadratic equalities make sure that $P'$ is similar to~$P$.

A global optimum of the optimization problem gives us a largest polytope $P'$ as desired. It might happen that there are combinatorially different optimal solutions to our problem. 
The goal in Section~\ref{loesung} is to identify \emph{one} of the optimal solutions. 
From that we can deduce the optimal dilation factor and hence answer the question: how large is the largest polytope $P'$ similar to $P$ and contained in~$Q$. 
We do not explain in full generality in what combinatorially different ways $P'$ can then be contained in $Q$, but rather describe one possible inclusion. 
\subsubsection{Improved formulation}
The above formulation for Problem~\ref{prob} is particularly simple and straightforward. 
However an equivalent formulation using less variables and less quadratic constraints can be obtained as follows.

Choose an affine basis from the set of vertices of~$P$. 
For the optimization problem we can then only take those variables $v_{ij}$, such that $w_i$ belong to that affine basis and substitute all occurrences of other variables by linear combination of the former. 
These linear combinations can be obtained from the vertices of $P$, using the fact that we chose an affine basis. 
Using this substitution, we have $(p+1)q+1$ variables in total and this number only depends on the dimensions of $P$ and $Q$ and not on the number of vertices of~$P$.

In order to obtain less quadratic constraints we also focus on the chosen affine basis: it is enough to make sure that all the distances between all pairs of two vectors in the affine basis are all scaled by the same factor $\sqrt{s}$. 
Since there are  $q+1$ vectors in the affine basis, we obtain $\binom{q+1}{2}=\frac{1}{2}(q+1)(q+2)$ quadratic equations. 
Counting the number of linear equations we see that there are $nm$ many, independent of the dimension of~$Q$. 

An axis aligned bounding box for $Q$ gives bounds on the variables $v_{ij}$. 
We can trivially include a copy of $P$, whose circumsphere coincides with the in sphere of $Q$, so a lower bound for $s$ would be the Keplerian ratio
\[s\geq\left(\frac{\text{circumradius of }Q}{\text{inradius of }P}\right).\]
In a similar way we could give an upper bound for $s$, but in view of the objective function this does not seem necessary.

The equations used in setting up Problem~\ref{prob} depend on the position of~$Q$. If many of the defining hyperplanes for $Q$ are parallel to many coordinate axes, then less variables are used in the linear equations. 
Also the choice of an affine basis of $P$ might influence the number of variables used in the equations.

The precision for the input of the polytopes should be higher than the desired precision, when solving Problem~\ref{prob} with a solver numerically.  

If $P$ and $Q$ possess symmetry one can use this symmetry to get additional constraints. 
For example if $P$ and $Q$ are centrally symmetric, then it suffices to search a maximal $P'$ among those copies of $P$ which are concentric with~$Q$. 
See \cite[\textsc{Observation} p.~288]{C80} for a simple proof.

If  $P$ and $Q$ are regular polytopes, one can say without loss of generality that one vertex of $P'$ must lie in one face of~$Q$.

\subsubsection{Solving the optimization problem numerically}
In order to solve Problem~\ref{prob} numerically we can use SCIP, which is a solver for mixed integer non-linear programming. 
This solver uses branch and bound techniques in order to find a \emph{global} optimum within a certain precision; see \cite{A09} and \cite{ABKW08} for details.
We don't use SCIP's capability to handle integer variables, since all of our variables are continuous.

\subsection{From numerical solutions to exact solutions}
\subsubsection{Setting up the quadratic system}

We obtain approximate results for the global optimum Problem \ref{prob}, with a certain precision
, let's call the resulting polytope~$\widetilde{P}$. 
The goal is to derive exact values for the coordinates of a polytope $P'$ which in indistinguishable from $\widetilde{P}$ in the approximation within the precision. 

We can identify the vertices of $\widetilde{P}$ that lie in a face of~$Q$. 
If $\widetilde{P}$ has been calculated with sufficiently high precision (see assumptions in Section~\ref{limits}) $\widetilde{P}$ will satisfy the same vertex-face incidences an optimal solution~$P'$.
In fact $P'$ is given by the real solution of a system of quadratic equations, which is derived from these incidences. An approximate real solution of this system is given by~$\widetilde{P}$.

\subsubsection{Solving the quadratic system}\label{method}
A numerical solution to this quadratic system with arbitrary precision can be obtained using Newton's method, and a solution $P'$ to Problem~\ref{prob} gives a good starting point.
If all the defining hyperplanes of $Q$ are defined in terms of algebraic numbers, solutions of the quadratic system must be algebraic. 
In case the system obtained in this way is to complicated to be solved by hand or automatically by a computer algebra system, we can attempt to find solutions by using the following three-step approach. 
We already have an approximate real solution given by~$\widetilde{P}$.

\begin{enumerate}
 \item[Step 1]\label{step1} Numerically approximate the solution to high precision, for example using multi-dimensional Newton's method
 \item[Step 2]\label{step2} For each variable guess the algebraic number close to the approximation using integer relation algorithms such as LLL (\cite{LLL82}).
 \item[Step 3]\label{step3} Verify the solution by exact calculation in the field of real algebraic numbers.
\end{enumerate}
We can expect to find solutions, if they are algebraic numbers with minimal polynomials of low degree and small coefficients. See Section~\ref{loesung} for two successful application of this method. 
This method can be in principle applied to any given system of equations with algebraic solutions, for which we can obtain high precision numerical approximate solutions.
\subsection{Limitations of the method}\label{limits}

The solver SCIP, which can be used for solving Problem~\ref{prob} finds a \emph{global} optimum, but the calculations are done only with a certain prescribed precision. 
In general it might be the case that exists a maximizer $P'$, which attains the maximal dilation factor $\sqrt{s}$ and a second locally maximal feasible solution $P''$, with dilation factor $\sqrt{s-\varepsilon}$, for a small~$\varepsilon>0$. 
Indeed it is possible to construct examples of $P$ and $Q$ where this is the case for arbitrarily small $\varepsilon$, take for example $P$ and $Q$ to both be the same rectangle with almost equal side length. 
Hence in order to make sure that we have indeed found an optimal solution to Problem~\ref{prob}, we make the following assumptions.
\begin{assumption}\label{assu1}
 The solution $\widetilde{P}$ to Problem~\ref{prob} has sufficient precision such that there is only one local maximum $P'$ near~$\widetilde{P}$.
\end{assumption}
\begin{assumption}\label{assu2}
 Problem~\ref{prob} has been solved with sufficient precision such that the dilation factor $\sqrt{s}$ of the local maximum $P'$ near $\widetilde{P}$ is the \emph{global} maximum. 
\end{assumption}
\begin{assumption}\label{assu3}
 Problem~\ref{prob} has been solved with sufficient precision such that $\widetilde{P}$ and the local maximum  $P'$ near $\widetilde{P}$ satisfy the same vertex-face incidences with~$Q$. 
\end{assumption}

The precision necessary for the solution to satisfy these properties depends on $P$ and $Q$ and since there exist examples where the global maximum and the second largest local maximum are arbitrarily close it is in general not possible to prescribe the precision necessary for Assumptions~\ref{assu1}-\ref{assu3} to hold.

Assumptions~\ref{assu1}-\ref{assu3} also deal with possible numerical mistakes or bugs of a solver for Problem~\ref{prob}.

If Assumptions~\ref{assu1} and \ref{assu2} hold and we can, because of Assumption~\ref{assu3} identify an exact algebraic solution near $P'$, this will be a maximizer of the problem. 
In any case, even if the assumptions do not hold, we get a lower bound if we can solve system derived from the approximate solution~$P'$.

In the calculations in Section~\ref{loesung} we do not attempt to prove that Assumptions~\ref{assu1}-\ref{assu3} hold, but we state the precision which was used to solve the problems.
In this sense our calculations below do not prove optimality but provide putatively optimal results.

\section{Results}\label{loesung}
\subsection{Inclusions of platonic solids}
When each of $P$ and $Q$ is taken to be one of the $5$ platonic solids, i.e.\ regular three-dimensional polyhedra, we can consider $20$ non-trivial inclusions. 
Croft found optimal pairs in $14$ out of these $20$ cases and proved optimality in \cite{C80}.
In the following we assume that the regular three-dimensional polyhedron $Q$ has side length~$1$. We abbreviate tetrahedron, cube, octahedron, dodecahedron and icosahedron by $T$,$C$,$O$,$D$ and $I$ respectively and denote the golden ratio by~$\phi$.

With the methods described above we are able to confirm all the known cases and answer all six unknown cases. 
The solver used was SCIP version 3.1.0 with a precision set to $10^{-10}$. With the improved formulation described above the calculations for all 20 inclusions took a few hours on a single core of a Xeon CPU running at 3 GHz, using less than 8GB of RAM. Some cases were solved in less than a second.
\begin{table}[H]
\centering
\begin{tabular}{|m{.16\textwidth}|m{.16\textwidth}|m{.16\textwidth}|m{.16\textwidth}|m{.16\textwidth}|}
 \hline
&\input{./CT.tikz}& \input{./OT.tikz}& \input{./DT.tikz}& \input{./IT.tikz}\\
\hline
\input{./TC.tikz}&& \input{./OC.tikz}& \input{./DC.tikz}& \input{./IC.tikz}\\
\hline 
\input{./TO.tikz}& \input{./CO.tikz}&& \input{./DO.tikz}& \input{./IO.tikz}\\
\hline 
\input{./TD.tikz}& \input{./CD.tikz}& \input{./OD.tikz}&& \input{./ID.tikz}\\
\hline 
\input{./TI.tikz}& \input{./CI.tikz}& \input{./OI.tikz}& \input{./DI.tikz}&\\
\hline
\end{tabular} 
\caption{Maximal platonic solids included in a platonic solid}\label{alltikz}
\end{table}

\pagebreak
The tables below give decimal approximations and symbolic values of the side length of a largest copy of $P$ inside $Q$, where $P$ and $Q$ range over the platonic solids. 
For completeness we restate the results of Croft, he gives a similar but incomplete table:  \cite[p.~295]{C80}. We correct three typos in his table, the corresponding cells are \emph{emphasized}; new results are marked with a $\star$star.
\begin{table}[H]
\noindent\resizebox{1 \textwidth}{!}{
\begin{tabular}{|c|c|c|c|c|c|}
 \hline\backslashbox{\scriptsize$Q$}{\scriptsize$P$}
&$T$&$C$&$O$&$D$&$I$\\
 \hline
 $T$& & 0.29590654 & 0.50000000 & $\star$ 0.16263158 & 0.27009076 \\
 \hline
$C$&1.4142136 &  & 1.0606602 & \emph{0.39428348} & 0.61803399 \\
 \hline
$O$&1.0000000 & 0.58578644 &  &$\star$ 0.31340182 & 0.54018151 \\
 \hline
$D$& \emph{2.2882456} & 1.6180340 &  \emph{1.8512296} &  & $\star$ 1.3090170\\
 \hline
$I$&$\star$ 1.3474429 & $\star$ 0.93874890 & 1.1810180 & $\star$ 0.58017873 &\\ 
\hline
\end{tabular}}
\end{table}
\begin{table}[H]
\extrarowsep=2mm
\resizebox{1 \textwidth}{!}{
 \begin{tabu}{|c|c|c|c|c|c|}
  \hline\backslashbox{\scriptsize$Q$}{\scriptsize$P$}
&$T$&$C$&$O$&$D$&$I$\\\hline
$T$& & $\frac{1}{1+\frac{2}{3}\sqrt{3}+\frac{1}{2}\sqrt{6}}$& $\frac{1}{2}$ & $\star  d$ & $\frac{1}{\phi^2\sqrt{2}}$ \\ 
 \hline
$C$& $\sqrt{2}$ & & $\frac{3}{4}\sqrt{2}$ & \scriptsize$\frac{1}{\sqrt{2}\phi^3}(1-\frac{1}{2}\sqrt{10}+\frac{1}{2}\sqrt{2}+\sqrt{5})$ & $\frac{1}{\phi}$ \\ 
 \hline
$O$& $1$ & $2-\sqrt{2}$ &&\large$\star\frac{(25\sqrt{2})-(9 \sqrt{10})}{22}$ & $\frac{\sqrt{2}}{\phi^2}$ \\ 
 \hline
$D$& $\phi\sqrt{2}$ &$\phi$&$\frac{\phi^2}{\sqrt{2}}$ & & \large$\star\frac{1}{2\phi}+1$ \\ 
 \hline
$I$&$\star  t$ &\large $\star\frac{5+7\sqrt{5}}{22}$ &\scriptsize$\frac{1}{2}(1-\frac{1}{2}\sqrt{10}+\frac{1}{2}\sqrt{2}+\sqrt{5})$&\large  $\star\frac{15-\sqrt{5}}{22}$&\\
 \hline
 \end{tabu}}\caption*{
\noindent$\phi = \text { golden ratio }$\\
  $t = \text{ zero near } 1.3 \text{ of } 5041x^{32} - 1318386 x^{30} + 60348584 x^{28} - 924552262 x^{26} + 5246771058 x^{24} - 15736320636   x^{22}  + 29448527368 x^{20} - 37805732980 x^{18}\\ + 35173457839 x^{16} - 24298372458 x^{14} + 12495147544 x^{12} - 4717349124x^{10}\\ + 1256858478 x^8- 217962112 x^6+21904868 x^4 -  1536272 x^2  + 160801$\\
  $d =  \text{ zero near } 0.16 \text{ of }  4096x^{16} - 3701760x^{14} + 809622720x^{12} - 17054118000x^{10} + 79233311025x^8 - 94166084250x^6 + 31024053000x^4 - 3236760000x^2 + 65610000$
}
\end{table}
  
\begin{figure}[H]
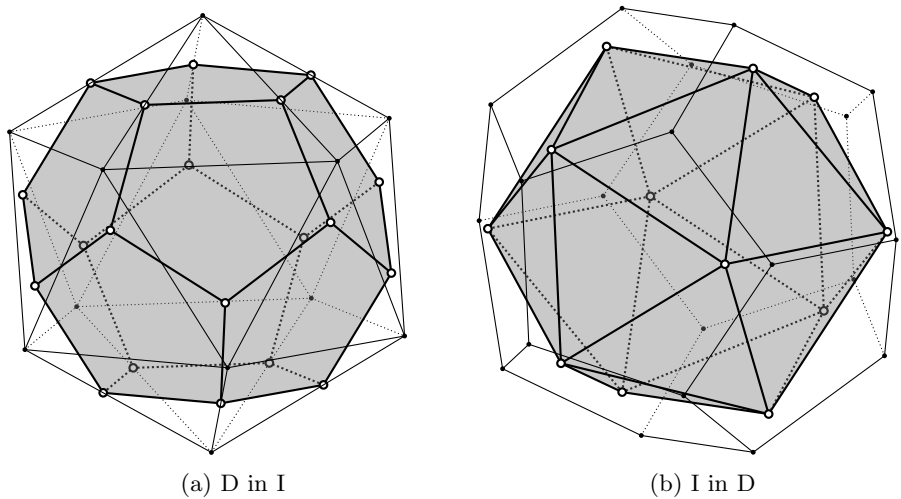

\centering
\begin{subfigure}[b]{0.49\textwidth}
\input{./DI2.tikz}
\subcaption{D in I}
\end{subfigure}~
\begin{subfigure}[b]{0.49\textwidth}
\input{./ID2.tikz}
\subcaption{I in D}
\end{subfigure}
\caption{Self reciprocal cases}
\end{figure}

\noindent For the $6$ previously unknown cases we give a description of an optimal position.

\subsubsection{Dodecahedron in icosahedron} 
For $D$ in $I$ we are in a concentric situation. The five vertices of one face of $D$ lie on the five edges of $I$ incident to a common vertex, one on each. 
The five vertices of the opposite face of that face of $D$ also lie on five edges of $I$ incident to a common vertex, namely the vertex of $I$ antipodal to the one mentioned before. The other ten vertices of $D$ lie in the interior of faces of~$I$. 
The side length is 
$$\frac{15-\sqrt{5}}{22}\approx 0.58017873.$$
  \subsubsection{Icosahedron in dodecahedron}
For $I$ in $D$ we are also in a concentric situation; each of the $12$ vertices of $I$ lies in the interior of one of the $12$ faces of $D$ and in each face of $D$ there is one vertex of~$I$. 
Let's position $D$ in the usual fashion such that $6$ of its edges are parallel to the $3$ coordinate axes. 
To each of the $12$ vertices on these edges of $D$ we associate the unique face which contains one but not the other vertex of the edge in its boundary. 
This gives us pairs $v,f$ of vertices and faces of~$D$. 
For each pair $v,f$ a vertex of $I$ lies on the bisector of $f$, which goes through $v$ and its position on the bisector is the point where the bisector is divided in two parts, such that the larger part has $\frac{\phi}{2}$ the length of the whole bisector.
The position of the vertex of $I$ is closer to $v$ and the absolute distance to $v$ is $(1-\frac{\phi}{2})\cdot \frac{1}{2}\sqrt[4]{5}\phi^{\frac{3}{2}}=\frac{\sqrt[4]{5}}{4\sqrt{\phi}}$. 
(Remember we assume that $D$ has side length $1$ which results in a bisector of length $\frac{1}{2}\sqrt[4]{5}\phi^{\frac{3}{2}}$.)
The edge length of $I$ obtained in this way is
\[\frac{1}{2\phi}+1\approx 1.3090170.\]
\begin{figure}[H]
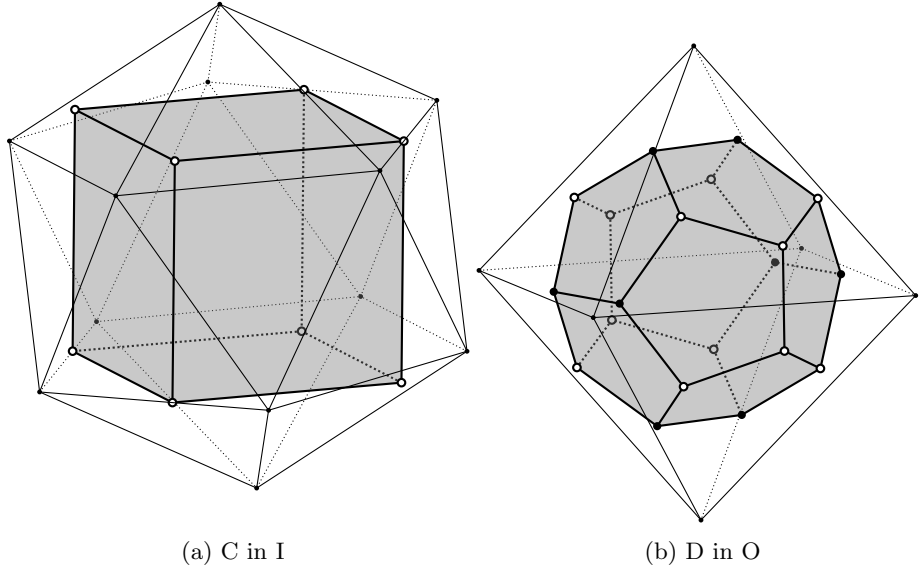

\centering
\begin{subfigure}[b]{0.49\textwidth}
\input{./CI2.tikz}
\subcaption{C in I}
\end{subfigure}~
\begin{subfigure}[b]{0.49\textwidth}
\input{./DO2.tikz}
\subcaption{D in O}\label{f2b}
\end{subfigure}
\caption{Two reciprocal cases}
\end{figure}
\subsubsection{Cube in icosahedron}
Also concentric. For $C$ in $I$, two vertices of one edge of $C$ lie in the interior of two adjacent edges in $I$, which are not contained in the same face. And the vertices of the antipodal edge of this edge in $C$ lie in the interior of the corresponding antipodal edges in $I$. The other 4 edges of $C$ lie in the interior of faces of $I$.
The side length is $$\frac{5+7\sqrt{5}}{22}\approx 0.93874890.$$
\subsubsection{Dodecahedron in octahedron}
 Again this is a concentric situation. Put two opposite edges of $D$ in a hyperplane spanned by $4$ vertices of $O$. 
 Four faces of $O$ each contain an edge of $D$ and the other four faces of $O$ each contain only one vertex of $D$. 
 The incidences can be seen in Figure~\ref{f2b}; vertices of $D$ which lie in the interior of a face of $O$ are marked white. 
 See the considerations about reciprocity below.
\noindent For $D$ in $O$ the maximum is
\[\frac{(25\sqrt{2})-(9 \sqrt{10})}{22}\approx0.31340182.\]
 \subsubsection*{Reciprocity of \texorpdfstring{$C\subset I$}{C in I} and \texorpdfstring{$D\subset O$}{D in O}} 
 If $P\subset Q$ are concentric and $P$ is maximal in $Q$ we can take polar reciprocals and get $Q^\circ\subset P^\circ$, such that $Q^\circ$ is maximal in $P^\circ$. 
 Since $C^\circ=O$ and $I^\circ=D$, we can check that the two previous cases are reciprocal: 
 \[\frac{(25\sqrt{2})-(9 \sqrt{10})}{22}\left(\frac{\phi^3}{\sqrt{2}}\right)=\frac{5+7\sqrt{5}}{22}.\]
 Concentric $C$ and $D$, which are reciprocals with respect to the unit sphere have the product of their edge lengths constant, namely $2\sqrt{2}$. 
 Similarly for concentric, reciprocal $I$ and $D$ this product equals $\frac{4}{\phi^3}$.
 The factor $\frac{\phi^3}{\sqrt{2}}$ is the quotient of these two numbers. 
\begin{figure}[H]
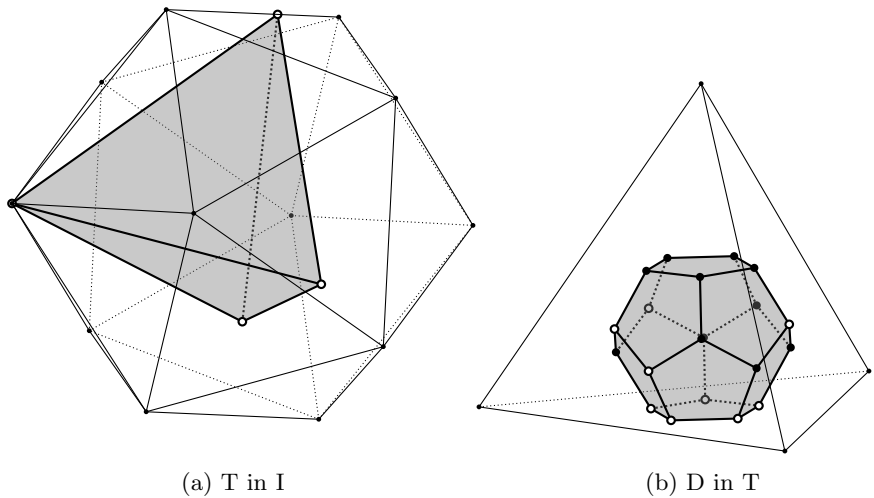

\centering
\begin{subfigure}[b]{0.49\textwidth}
\input{./TI2.tikz}
\subcaption{T in I}\label{f3a}
\end{subfigure}~
\begin{subfigure}[b]{0.49\textwidth}
\input{./DT2.tikz}
\subcaption{D in T}\label{f3b}
\end{subfigure}
\caption{Two cases with more involved solutions}\label{f3}
\end{figure}

\subsubsection{Tetrahedron in icosahedron}
The incidences of the $T$ in $I$ are best seen in Figure~\ref{f3a}: one vertex of $T$ coincides with one vertex $v$ of $I$, another vertex of $T$ lies on an edge of $I$, which is neither incident to the vertex $v$ nor its antipode, and the two remaining vertices lie in the interior of faces of $I$.

While in this case the resulting system can be somewhat automatically solved by the computer algebra system Mathematica 9 (while version 8 was not able to perform the calculation), we use the methods described in Section~\ref{method}. 
We choose two variables each for the barycentric coordinates for the two vertices in the interior of faces of $I$ and one variable for barycentric coordinates for the vertex in the interior of an edge of $I$. 
Together with a variable $t$ for the side length of $T$, i.e.\ the dilation factor, this results in a system of $6$ quadratic equations in $6$ variables. The $6$ equations confirm that all $6$ edges are of length $t$. 
We use the open source computer algebra system \emph{sage}, \cite{sage}. For the newton method, i.e.\ Step~1 we use scipy, \cite{scipy}, and for the integer relation, i.e.\ Step 2 PARI, \cite{PARI2} is used.  It is sufficient to obtain 800 decimal digits in Step~1 of the method described in Section~\ref{method} in order to obtain the exact values for the variables in Step~2.
The exact edge length is the zero near $1.3474429$ of this polynomial:
\[  \begin{array}{l}5041 t^{32}-1318386 t^{30}+60348584 t^{28} -924552262 t^{26}+5246771058 t^{24}\\
		    -15736320636   t^{22}  +29448527368 t^{20}-37805732980 t^{18}+35173457839 t^{16} \\
 		-24298372458 t^{14}+12495147544 t^{12}-4717349124t^{10}	+1256858478 t^8\\-217962112 t^6
		+21904868 t^4-1536272 t^2   +160801.\end{array}\]

\subsubsection{Dodecahedron in tetrahedron}
The incidences are best seen in Figure~\ref{f3b}: a complete face of $D$ is contained in one face of $T$, two vertices of $D$ lie in another face of $T$ and the two other faces of $T$ contain one vertex of $D$ each. 
We choose a variable $d$ for the side length of $D$ and four additional variables that describe the position of the vertices of $D$ that lie in a face of $T$, which is not the face that contains a complete face of $D$.  
Making sure that the edges between these four vertices have the correct length results again in a system of $6$ quadratic equations with $5$ variables, which can be successfully solved as in the previous case.
In this case 350 decimal digits suffice to find solutions in the field of real algebraic numbers. 
The exact edge length is the zero near $0.16263158$ of this polynomial:
\[ \begin{array}{l}  4096d^{16} - 3701760d^{14} + 809622720d^{12} - 17054118000d^{10} + 79233311025d^8 - \\94166084250d^6 + 31024053000d^4 - 3236760000d^2 + 65610000.
      \end{array}\]

\section{Further applications}
Possibly interesting situations where the method of this paper could be applied include the following cases.
\begin{enumerate}[a)]
 \item\label{e1} Take $P$ and $Q$ to be (regular) polygons.
 \item Take $P$ and $Q$ to be regular polytopes of dimension greater than~$3$.
 \item Take $P$ to be a $n$-cube and $Q$ an $m$-cube with $n<m$.
 \item Take $P$ to be a regular $n$-simplex and $Q$ an $m$-cube with $n\leq m$.
 \item Take $P$ to be a regular $n$-simplex and $Q$ an regular $m$-simplex with $n<m$.
 \item Take $Q$ to be any polytope and $P$ some projection of $Q$.
\end{enumerate}
For the first case, i.e.\ finding the largest regular $n$-gon in a regular $m$-gon, the author has checked the conjecture of Dilworth and Mane \cite[Section 9]{DM10} for coprime $m$ and $n$ up to a precision of $10^{-10}$ for all pairs $m,n$ with $m,n\leq120$.

It is possible to modify Problem~\ref{prob} in order to solve similar packing problems. 

\section*{Acknowledgements}\label{ackref}
  I would like to thank Ambros Gleixner, Günter M. Ziegler, Hartmut Monien, Louis Theran and Peter Bürgisser for fruitful discussions. 
\bibliographystyle{alpha}

\bibliography{lit.bib}

\end{document}